\documentclass[12pt]{article}

\usepackage{amsmath,amssymb,amsthm,cite,enumitem,mathtools}
\usepackage[top=3cm, bottom=3cm, left=2.8cm, right=2.8cm]{geometry}

\usepackage[osf,sups,scaled=.97]{XCharter}
\usepackage[libertine,bigdelims,vvarbb,scaled=1.03]{newtxmath}
\usepackage[cal=euler,frak=euler]{mathalfa} % mathcal

\title{Faà di Bruno Hopf algebras}

\author{Héctor Figueroa,$^1$
José M. Gracia-Bondía$^{2,3}$
and Joseph C. Várilly$^1$% 
\footnote{Email: joseph.varilly@ucr.ac.cr}
\\[6pt]
{\footnotesize $^1$ Escuela de Matemática, 
Universidad de Costa Rica, San José 11501, Costa Rica}
\\[3pt]
{\footnotesize $^2$ CAPA and Departamento de Física Teórica,
Universidad de Zaragoza, Zaragoza 50009, Spain}
\\[3pt]
{\footnotesize $^3$ Laboratorio de Física Teórica y Computacional,
Universidad de Costa Rica, San Pedro 11501, Costa Rica}
}

\date{}

\DeclareMathOperator{\End}{End}     %% endomorphism space
\DeclareMathOperator{\Hom}{Hom}     %% homomorphism space
\DeclareMathOperator{\sh}{sh}       %% shuffles

\newcommand{\al}{\alpha}            %% short for \alpha
\newcommand{\Dl}{\Delta}            %% coproduct
\newcommand{\dl}{\delta}            %% short for \delta
\newcommand{\Ga}{\Gamma}            %% short for \Gamma
\newcommand{\la}{\lambda}           %% short for \lambda
\newcommand{\sg}{\sigma}            %% short for  \sigma

\newcommand{\bN}{{\mathbb{N}}}      %% natural numbers
\newcommand{\bR}{{\mathbb{R}}}      %% real numbers

\newcommand{\sA}{\mathcal{A}}       %% Connes-Moscovici Lie algebra
\newcommand{\sF}{\mathcal{F}}       %% Faà di~Bruno Hopf algebra
\newcommand{\sP}{\mathcal{P}}       %% hereditary family
\newcommand{\sU}{\mathcal{U}}       %% enveloping algebra

\newcommand{\alg}{{\mathrm{alg}}}   %% algebraic version
\newcommand{\ivth}{{(\mathrm{iv})}} %% fourth derivative

\newcommand{\del}{\partial}         %% short for \partial
\newcommand{\otto}{\leftrightarrow} %% bijection
\newcommand{\ox}{\otimes}           %% tensor product
\newcommand{\stroke}{\mathbin|}     %% restricted patition
\newcommand{\tPi}{\widetilde{\Pi}}  %% type of partition
\newcommand{\wt}{\widetilde}        %% short for  \widetilde
\newcommand{\xyx}{\times\cdots\times} %% repeated Cartesian product
\renewcommand{\:}{\colon}           %% colon in  f: A -> B

\newcommand{\set}[1]{\{\,#1\,\}}    %% set notation
\newcommand{\word}[1]{\quad\text{#1}\quad} %% well-spaced word(s)

\newcommand{\duo}[2]{\langle#1,#2\rangle} %% dual pairing
\newcommand{\stirlingtwo}[2]{\genfrac\{\}{0pt}{}{#1}{#2}} % n subset k

\makeatletter
\renewcommand{\section}{\@startsection{section}{1}{\z@}%
							 {-3.25ex \@plus -1ex \@minus -.2ex}%
							 {1.5ex \@plus.2ex}%
							 {\normalfont\large\bfseries}}
\makeatother

\begin{document}

\maketitle

\begin{abstract}
This is a short review on the Faà di Bruno formulas, implementing
composition of real-analytic functions, and a Hopf algebra associated
to such formulas. This structure allows, among several other things, a
short proof of the Lie--Scheffers theorem, and relating the Lagrange
inversion formulas with antipodes. It is also the maximal commutative
Hopf subalgebra of the one used by Connes and Moscovici to study
diffeomorphisms in a noncommutative geometry setting. The link of Faà
di~Bruno formulas with the theory of set partitions is developed in
some detail.
\end{abstract}

% \keywords{Faà di Bruno formula, Hopf algebras, partitions}
% \subjclass{16T30, 05A17}

\section{The Faà di Bruno formula} % 1
\label{sec:FdB-formula}

Faà di~Bruno (Hopf, bi)algebras appear in several branches of
mathematics and physics, and may be introduced in several ways. Here
we start from the group~$G$ of formal exponential power series like
$$
f(t) = \sum_{n=1}^\infty \frac{f_n}{n!} \,t^n,
$$
with $f_1 > 0$. (In view of Borel's lemma, one may regard them as local
representatives of orientation-preserving diffeomorphisms of~$\bR$
leaving $0$~fixed. The question of analyticity is considered below.)

On this group of power series we consider the coordinate functions
$$
a_n(f) := f_n = f^{(n)}(0), \quad n \geq 1.
$$
We wish to compute $h_n = a_n(h)$, where $h$ is the composition
$f \circ g$ of two such diffeomorphisms, in terms of the $f_n$
and~$g_n$. Now,
$$
h(t) = \sum_{k=1}^\infty \frac{f_k}{k!}
\biggl( \sum_{l=1}^\infty \frac{g_l}{l!} \,t^l \biggr)^k.
$$
To compute the $n$th coefficient $h_n$ we need only consider the sum
up to $k = n$, since the remaining terms contain powers of~$t$ higher
than~$n$. From Cauchy's product formula,
$$
h_n = \sum_{k=1}^n \frac{f_k}{k!} \sum_{l_i\geq 1,\,l_1+\cdots+l_k=n}
\frac{n!\,g_{l_1}\cdots g_{l_k}}{l_1!\cdots l_k!}.
$$
If among the $l_i$ there are $\la_1$ copies of~$1$, $\la_2$ copies
of~$2$, and so on, then the sum $l_1 +\cdots+ l_k = n$ can be
rewritten as
\begin{equation}
\la_1 + 2\la_2 +\cdots+ n\la_n = n,  \word{with}  
\la_1 +\cdots+ \la_n = k.
\label{eq:la-count} % (1)
\end{equation}
Since there are $k!/\la_1!\cdots\la_n!$ contributions from~$g$ of this
type, it follows that
\begin{equation}
h_n = \sum_{k=1}^n f_k \sum_\la \frac{n!}{\la_1! \cdots \la_n!}
\frac{g_1^{\la_1} \cdots g_n^{\la_n}}
     {(1!)^{\la_1} \,(2!)^{\la_2} \cdots (n!)^{\la_n}}
=: \sum_{k=1}^n f_k \, B_{n,k}(g_1,\dots,g_{n+1-k}),
\label{eq:wonder-series} % (2)
\end{equation}
where the sum $\sum_\la$ runs over the sequences
$\la = (\la_1,\dots,\la_n) \in \bN^\bN$
satisfying~\eqref{eq:la-count}, and the $B_{n,k}$ are called the
(partial, exponential) \textit{Bell polynomials}. Usually these are
introduced by the expansion
$$
\exp\biggl( u \sum_{m\geq 1} x_m \frac{t^m}{m!} \biggr)
= 1 + \sum_{n\geq 1} \frac{t^n}{n!}
\biggl[ \sum_{k=1}^n u^k B_{n,k}(x_1,\dots,x_{n+1-k}) \biggr],
$$
which is a particular case of~\eqref{eq:wonder-series}. Each $B_{n,k}$
is a homogeneous polynomial of degree~$k$. (This is a good moment to
declare that the scalar field $\bR$ may be replaced by any commutative
field of characteristic zero.)

Formula~\eqref{eq:wonder-series} can be recast as
\begin{equation}
h^{(n)}(t) = \sum_{k=1}^n \sum_\la \frac{n!}{\la_1!\cdots\la_n!}
f^{(k)}(g(t)) \biggl( \frac{g^{(1)}(t)}{1!} \biggr)^{\la_1}
\biggl( \frac{g^{(2)}(t)}{2!} \biggr)^{\la_2} \! \cdots
\biggl( \frac{g^{(n)}(t)}{n!} \biggr)^{\la_n}.
\label{eq:remix} % (3)
\end{equation}
Expression~\eqref{eq:remix} is the famous formula attributed to Faà
di~Bruno (1855, 1857), who in fact followed previous authors; his
original contribution was a determinant form of it. Apparently
\eqref{eq:remix} goes back to Arbogast (1800); we refer the reader
to~\cite{Craik05} ---and references therein--- for these historical
matters. Note that if $g,f$ are differentiable up to order~$n$, then
$h$ is also differentiable up to order~$n$, and the expressions for
its derivatives hold.

The formula shows that the composition of two real-analytic functions
is real-analytic. Indeed, by use of \eqref{eq:wonder-series} or
\eqref{eq:remix} with $f(t) = \sum_{k=1}^\infty t^k = t/(1 - t)$ and
$g(t) = \sum_{l=1}^\infty x t^l = xt/(1 - t)$ one sees that
\begin{equation}
\sum_{k=1}^n \sum_\la \frac{k!}{\la_1!\la_2!\cdots\la_n!}\, x^k
= x(1 + x)^{n-1}.
\label{eq:keep-calm} % (4)
\end{equation}
Now, a smooth function $g$ on an open interval $I \subseteq \bR$ is
analytic~\cite[Chap.~1]{KrantzP02} if and only if for each $y \in I$
there is an open interval $J_y$ with $y \in J_y \subseteq I$ and
constants $A$, $B$ such that
\begin{equation}
|g^{(j)}(t)| \leq \frac{A\,j!}{B^j}  \word{for all} t \in J_y
\label{eq:tailored} % (5)
\end{equation}
which guarantees local uniform convergence of the Taylor series
of~$g$. Assume further that $g$~takes values in an open interval on
which the smooth function $f$ is defined. If $f$~is also analytic with
$|f^{(m)}(s)| \leq C\,m!/D^m$ for all~$m$ at $s = g(t)$, it follows
from~\eqref{eq:remix} and \eqref{eq:keep-calm} that
$$
|h^{(n)}(t)|
\leq \sum_{k=1}^n \sum_\la \frac{n!}{\la_1!\cdots\la_n!}
\frac{C\,k!}{D^k} \biggl(\frac{A}{B}\biggr)^{\la_1} \! \cdots
\biggl(\frac{A}{B^n}\biggr)^{\la_n} \!
= n!\,\frac{C}{B^n} \frac{A}{D} \biggl(1 + \frac{A}{D}\biggr)^{n-1} \!
= \frac{E\,n!}{F^n},
$$
with $E = AC/(A + D)$ and $F = BD/(A + D)$. Hence $f \circ g$ is
analytic on the domain of~$g$.

\section{Hopf algebras} % 2
\label{sec:Hopf-algebras}

Introduce the notation, with \eqref{eq:la-count} understood:
$$
\binom{n}{\la;k} := \frac{n!}{\la_1!\cdots\la_n!\,
(1!)^{\la_1}(2!)^{\la_2}\dots(n!)^{\la_n}}.
$$
A \textit{Hopf algebra} dual to~$G$ is obtained when we define a
\textit{coproduct}~$\Dl$ on the polynomial algebra 
$\bR[a_1,a_2,\dots]$ of coordinate functions
by requiring that $\Dl a_n(g,f) = a_n(f \circ g)$, or equivalently,
$a_n(f \circ g) = m(\Dl a_n(g \ox f))$ where $m$ means multiplication.
This entails that
$$
\Dl a_n = \sum_{k=1}^n \sum_\la \binom{n}{\la;k}
a_1^{\la_1} a_2^{\la_2} \dots a_n^{\la_n} \ox a_k.
$$
The unnecessary flip of~$f$ and~$g$ is traditional. This \textit{Faà
di~Bruno bialgebra}, so called by Joni and Rota~\cite{JoniR82}, is
commutative but not cocommutative. Since $a_1$ is a grouplike element,
it must be invertible for this to be a Hopf algebra. For that, one
must either adjoin an inverse $a_1^{-1}$, or else put $a_1 = 1$, as we
do from now on. That is, we consider only the subgroup $G_1$ of
diffeomorphisms tangent to the identity at~$0$. The first instances of
the coproduct are, accordingly,
\begin{align}
\Dl a_2 &= a_2 \ox 1 + 1 \ox a_2,
\notag \\
\Dl a_3 &= a_3 \ox 1 +  1 \ox a_3 + 3a_2 \ox a_2,
\notag \\
\Dl a_4 &= a_4 \ox 1 +  1 \ox a_4 + 6a_2 \ox a_3
+ (3a_2^2 + 4a_3) \ox a_2,
\notag \\
\Dl a_5 &= a_5 \ox 1  + 1 \ox a_5 + 10a_2 \ox a_4
+ (10a_3 + 15a_2^2) \ox a_3 + (5a_4 + 10a_2a_3) \ox a_2.
\label{eq:delta-bfb} % (6)
\end{align}
The resulting graded connected Hopf algebra $\sF$ is called the
\textit{Faà di~Bruno Hopf algebra}; the degree~$\#$ being given by
$\# a_n = n - 1$.

Consider the graded dual Hopf algebra~$\sF'$. Its space of primitive
elements has a basis $\{a'_n : n \geq 2\}$ defined by
$\duo{a'_n}{a_m} = \dl_{nm}$ and
$\duo{a'_n}{a_{m_1}a_{m_2}\dots a_{m_r}} = 0$ for $r > 1$. Their
product is given by the duality recipe
$\duo{b'c'}{a} := \duo{b'\ox c'}{\Dl a}$, leading to:
$$
a'_n a'_m = \binom{m-1+n}{n} a'_{n+m-1} + (1 + \dl_{nm})(a_n a_m)'.
$$
In particular, taking $b'_n := (n+1)!\,a'_{n+1}$ for $n \geq 1$, we
are left with the commutator relations
\begin{equation}
[b'_n, b'_m] = (m - n) b'_{n+m}.
\label{eq:so-nice} % (7)
\end{equation}
The Milnor--Moore theorem implies that $\sF'$ is isomorphic to the
enveloping algebra of the Lie algebra $\sA$ spanned by the~$b'_n$ with
these commutators.

A curious consequence of~\eqref{eq:so-nice} is that the space $P(\sF)$
of primitive elements of~$\sF$ just has dimension~$2$. Indeed, 
$P(\sF) = (\bR 1 \oplus {\sF'_+}^2)^\perp$, where $\sF'_+$ is the
augmentation ideal of~$\sF'$. But \eqref{eq:so-nice} entails that
there is a basis of~$\sF'$ made of products, except for its first two
elements: therefore, $\dim P(\sF) = 2$. A basis of $P(\sF)$ is given
by $\{a_2, a_3 - \frac{3}{2} a_2^2\}$. The second of these corresponds
to the Schwarzian derivative, which is known~\cite{Hille76} to be
invariant under the projective group $PSL(2,\bR)$. Nonexistence of
more primitive elements of~$\sF$ is related to the affine linear and
Riccati equations being the only Lie--Scheffers
systems~\cite{LieS1893,Kalliope} over the real~line.

\medskip

The Faà di~Bruno Hopf algebra $\sF$ reappears as the maximal
commutative Hopf subalgebra of the (noncommutative geometry) Hopf
algebra $H$ of Connes and Moscovici~\cite{ConnesM98}. Their
description of~$\sF$ uses a different set of coordinates
$\dl_n(f) := [\log f'(t)]^{(n)}(0)$, $n \geq 1$. Since
$$
h(t) := \sum_{n\geq 1} \dl_n(f) \frac{t^n}{n!} = \log f'(t)
= \log\biggl( 1 + \sum_{n\geq 1} a_{n+1}(f)\,\frac{t^n}{n!} \biggr),
$$
it follows from formula~\eqref{eq:wonder-series}, for logarithm and
exponential functions respectively, that
\begin{align*}
\dl_n &= \sum_{k=1}^n (-1)^{k-1}(k-1)!\,B_{n,k}(a_2,\dots,a_{n+2-k})
=: L_n(a_2,\dots,a_{n+1}), 
\\
\shortintertext{inverted by} 
a_{n+1} &= \sum_{k=1}^n B_{n,k}(\dl_1,\dots,\dl_{n+1-k})
=: Y_n(\dl_1,\dots,\dl_n),
\end{align*}
where the $L_n$ and the $Y_n$ are respectively called the
\textit{logarithmic polynomials} and the (complete, exponential)
\textit{Bell polynomials}. In this way we get $\dl_1 = a_2$;
$\dl_2 = a_3 - a_2^2$; $\dl_3 = a_4 - 3a_2a_3 + 2a_2^3$; 
$\dl_4 = a_5 - 3a_3^2 - 4a_2a_4 + 12a_2^2a_3 - 6a_2^4$; and so on.
Since the coproduct is an algebra morphism, by use
of~\eqref{eq:delta-bfb} we may obtain the coproduct in the
Connes--Moscovici coordinates. For instance,
\begin{equation*}
\Dl\dl_4 = \dl_4 \ox 1 + 1 \ox \dl_4 + 6\dl_1 \ox \dl_3 + (7\dl_1^2 +
4\dl_2) \ox \dl_2 + (3\dl_1 \dl_2 + \dl_1^3 + \dl_3) \ox \dl_1.
\end{equation*}

It is not easy to find a closed formula for $\Dl(\dl_n)$ directly
from~\eqref{eq:delta-bfb}. Fortunately, through~$\sF'$ another method
is available. Using $B_{n,1}(a_2,\dots,a_{n+1}) = a_{n+1}$, one finds
that $\duo{b'_n}{\dl_m} = (n+1)!\,\dl_{n,m}$. Let~$A$ be the graded
free Lie algebra generated by primitive elements $X_n$, $n \geq 1$.
Its enveloping algebra $\sU(A)$ is the \textit{concatenation Hopf
algebra}. A linear basis for $\sU(A)$, indexed by all vectors with
positive integer components $\bar n = (n_1,\dots,n_r)$, is made of
products $X_{\bar n} := X_{n_1}X_{n_2}\dots X_{n_r}$, together with
the unit element $X_\emptyset = 1$. Its coproduct is
\begin{equation*}
\Dl(X_{\bar n}) := \sum_{\bar n^1,\bar n^2}
\sh_{\bar n}^{\bar n^1,\bar n^2} X_{\bar n^1} \ox X_{\bar n^2},
\end{equation*}
with $\sh_{\bar n}^{\bar n^1,\bar n^2}$ denoting the number of
shuffles of the vectors $\bar n^1$, $\bar n^2$ that 
produce~$\bar n$. Let
$u^{\bar n}$ denote a dual basis to $X_{\bar n}$; the graded dual of
$\sU(A)$ is the \textit{shuffle Hopf algebra} $H$ with product and
coproduct respectively given by
$$
u^{\bar n^1} u^{\bar n^2}
:= \sum_{\bar n} \sh_{\bar n}^{\bar n^1,\bar n^2} u^{\bar n}, \qquad
\Dl(u^{\bar n})
:= \sum_{\bar n^1\bar n^2=\bar n} u^{\bar n^1} \ox u^{\bar n^2},
$$   
where $\bar n^1\bar n^2$ is the concatenation of the vectors
$\bar n^1$, $\bar n^2$. The surjective morphism $\rho\: A \to \sA$
defined by $\rho(X_n) := b'_n$ extends, by the universal property of
enveloping algebras, to a surjective morphism 
$\rho\: \sU(A) \to \sF'$, whose transpose is the injective Hopf map
$\rho^t\: \sF \to H$ given by 
$\dl_n \mapsto \Ga_n := \dl_n \circ \rho$. We may thus regard $\sF$ as
a Hopf subalgebra of~$H$, and thereby compute the coproduct of~$\sF$
from that of~$H$. The argument may look circular, since we seem to
need an expression for the~$\Ga_n$, which in turn requires computing
$\Dl(\dl_n)$. But we can write
\begin{equation}
\duo{\Ga_m}{X_{\bar n}} = \duo{\dl_m}{\rho(X_{\bar n})}
= \duo{\dl_m}{b'_{n_1} \cdots b'_{n_r}}
= \duo{\Dl(\dl_m)}{b'_{n_1} \ox b'_{n_2} \cdots b'_{n_r}}.
\label{eq:what-is-Gamma} % (8)
\end{equation}
Thus, to compute $\Ga_n$, the only terms we need in the expansion of
$\Dl(\dl_n)$ are $\dl_n \ox 1 + 1 \ox \dl_n$ and the bilinear terms,
namely multiples of $\dl_i \ox \dl_j$; the remaining terms are of the
form 
$(\mathrm{constant})\,\dl_{i_1}^{r_1}\cdots \dl_{i_k}^{r_k}\ox \dl_j$,
where $r_1 i_1 +\cdots+ r_k i_k + j = n$. The bilinear part $B(\dl_n)$
may be computed by induction~\cite{ConnesM98} to be
\begin{equation}
B(\dl_n) = \sum_{i=1}^{n-1} \binom{n}{i-1} \dl_{n-i} \ox \dl_i.
\label{eq:bilinear-part} % (9)
\end{equation}
Substituting \eqref{eq:bilinear-part} repeatedly
in~\eqref{eq:what-is-Gamma}, one obtains
$$
\Ga_n = n! \sum_{\bar n : n_1+\cdots+n_r = n} C^{\bar n}\, u^{\bar n}
\word{with coefficients} 
C^{\bar n} := (n_r + 1) \prod_{i=2}^r (n_i +\cdots+ n_r).
$$
For instance, $\Ga_1 = 2u^1$ and 
$\Ga_3 = 12(2u^3 + u^{(2,1)} + 3u^{(1,2)} + 2u^{(1,1,1)})$.

Another calculation of the~$\Ga_n$ was sketched in~\cite{Menous05}; it
eventually allows to improve~\eqref{eq:bilinear-part} to
$$
\Dl(\dl_n) = \dl_n \ox 1 + 1 \ox \dl_n
+ \sum_{\bar n\in N_n} \frac{n!}{n_1!\dots n_r!}
\,K_{n_r}^{n_1,\dots,n_{r-1}}
\,\dl_{n_1} \dots \dl_{n_{r-1}} \ox \dl_{n_r},
$$
where $N_n := \set{\bar n : n_1+\cdots+n_r = n, \ r > 1}$ and,
mindful that $\binom{n_r}{k} = 0$ when $n_r < k$,
$$
K_{n_r}^{n_1,\dots,n_{r-1}}
= \sum_{k=1}^{r-1} \binom{n_r}{k} 
\sum_{\bar n^1\cdots\bar n^k = (n_1,\dots,n_{r-1})}
\frac{1}{r^1!\dots r^k!} 
\prod_{i=1}^k \frac{1}{1 + n_1^i +\cdots+ n_{r^i}^i}\,.
$$
For $r = 2$, this becomes $K^{n-i}_i = \frac{i}{1+n-i}$, thus the
coefficient of $\dl_{n-i} \ox \dl_i$ is
$\binom{n}{i} \frac{i}{1+n-i} = \binom{n}{i-1}$, as
in~\eqref{eq:bilinear-part}.

\medskip

From the combinatorial viewpoint, the Faà di~Bruno Hopf algebra is the
\textit{incidence Hopf algebra} corresponding to intervals formed by
partitions of finite sets. This is no surprise, since the coefficients
of a Bell polynomial $B_{n,k}$ just count the number of partitions of
$\{1,\dots,n\}$ into $k$ blocks. A \textit{partition} 
$\pi \in \Pi(S)$, of a finite set $S$ with $n$ elements, is a
collection $\{B_1,B_2,\dots,B_k\}$ of nonempty disjoint subsets,
called \textit{blocks}, such that $\bigcup_{i=1}^k B_i = S$. We simply
write $\pi \vdash n$ for such, with $|\pi|$ being the number of blocks
in~$\pi$.

We say that $\pi$ is of \textit{type} $(\al_1,\dots,\al_n)$ if exactly
$\al_i$ of these $B_j$ have $i$~elements; thus
$\al_1 + 2\al_2 +\cdots+ n\al_n = n$ and $\al_1 +\cdots+ \al_n = k$
\cite{Stanley78}. We say that $\pi$ \textit{refines} $\tau$, and write
$\{A_1,\dots,A_n\} = \pi \leq \tau = \{B_1,\dots,B_m\}$, if each~$A_i$
is contained in some~$B_j$. A subinterval
$[\pi,\tau] = \{\sg : \pi \leq \sg \leq \tau\}$ of the lattice $\sP$
of partitions of finite sets is isomorphic to the poset
$\Pi_1^{\la_1} \xyx \Pi_n^{\la_n}$, where 
$\Pi_j := \Pi(\{1,\dots,j\})$ and $\la_i$ blocks of $\tau$ are unions
of exactly $i$ blocks of~$\pi$. One assigns to each interval the
sequence $\la = (\la_1,\dots\!,\la_n)$ and declares two intervals
in~$\sP$ to be equivalent when their vectors $\la$ are equal. From the
matching $\wt{[\pi,\tau]} \otto \la \otto
\tPi_1^{\la_1} \tPi_2^{\la_2} \cdots \tPi_n^{\la_n}$ of equivalence
classes, one may regard the family $\wt\sP$ of equivalence classes as
the algebra of polynomials of infinitely many variables
$\bR[\tPi_1,\tPi_2,\dots]$. By means of the general theory of
coproducts for incidence bialgebras~\cite{Quaoar} one then recovers
the Faà di~Bruno algebra under the identifications $a_n \otto \tPi_n$.
The cardinality in the sense of category theory~\cite{BaezD01} of the
groupoid of finite sets equipped with a partition is given by
$$
\sum_{n=0}^\infty \sum_{k=1}^n \frac{1}{n!}\,B_{n,k}(1,\dots,1)
= e^{e-1}.
$$

The \textit{characters} of~$\sF$ form a group $\Hom_\alg(\sF,\bR)$
under the convolution operation of Hopf algebra theory. The action of
a character $f$ is determined by its values on the~$a_n$. The map
$f \mapsto f(t) = \sum_{n=1}^\infty f_n t^n/n!$, where
$f_n := \duo{f}{a_n}$, matches characters with exponential power
series over~$\bR$ such that $f_1 = 1$. This correspondence is an
\textit{anti-isomorphism} of groups: indeed, the convolution $f * g$
of $f,g \in \Hom_\alg(\sF,\bR)$ is given by
$$
\duo{f * g}{a_n} := m(f \ox g)\Dl a_n = \duo{g \circ f}{a_n}.
$$
This is just the $n$th coefficient of $h(t) = g(f(t))$. Also, the
algebra endomorphisms $\End_\alg(\sF)$ form a group under the
convolution of the unital algebra $\End(\sF)$ of linear endomorphisms.

The inverse under functional composition of an exponential series is
given by the reversion formula of Lagrange~\cite{Lagrange1770}, one of
whose forms~\cite{Comtet74} states that if $f$ and~$g$ are two such
series and if $f_1 = 1$, $f \circ g(t) = g \circ f(t) = t$, then
\begin{equation}
g_n = \sum_{k=1}^{n-1} (-1)^k B_{n-1+k,k}(0,f_2,f_3,\dots).
\label{eq:no-cancelation} % (10)
\end{equation}
Now, the inverse under convolution of $f \in \End_\alg(\sF)$ is
$g = f\circ S$, with $S$ the antipode map of~$\sF$. The
multiplicativity of~$f$ forces
$$
S(a_n) = \sum_{k=1}^{n-1} (-1)^k B_{n-1+k,k}(0,a_2,a_3,\dots).
$$
One may reverse the roles and prove the combinatorial
identity~\eqref{eq:no-cancelation} from Hopf algebra
theory~\cite{HaimanS89}. The real-analytic inverse function theorem is
stated in a completely similar way to the standard inverse function
theorem for differentiable functions, and it can be proved by use of
the formula of Faà di~Bruno, with the help of the
estimates~\eqref{eq:tailored}.

Use of partitions with special properties may lead to other incidence
algebras: for instance, if we restrict to noncrossing partitions, we
obtain a cocommutative Hopf algebra, with the commutative group
operation on characters essentially corresponding to Lagrange
reversion of the Cauchy product of reverted series~\cite{Stanley99}.

\section{Faà di Bruno formulas in several variables} % 3
\label{sec:ex-uno-plures}

To go to higher Faà di~Bruno formulas means to consider exponential
$N'$-series in $N'$ variables (``colours'') of the form
\begin{equation}
f(t_1,\dots,t_{N'}) = \biggl(
t_1 + \sum_{|\bar m|>1} f^1_{\bar m} \frac{t^{\bar m}}{\bar m!}\,,
t_2 + \sum_{|\bar n|>1} f^2_{\bar n} \frac{t^{\bar n}}{\bar n!}\,,
\dots,
t_{N'} + \sum_{|\bar p|>1} f^{N'}_{\bar p} \frac{t^{\bar p}}{\bar p!}
\biggr),
\label{eq:N'-series} % (11)
\end{equation}
where $\bar m,\bar n,\dots,\bar p \in \bN^{N'}$. The simplest way to
go about this is to rewrite Eq.~\eqref{eq:remix} in terms of
partitions:
\begin{equation}
(f\circ g)^{(n)}(t) = \sum_{\pi\vdash n} f^{(|\pi|)}(g(t))
\prod_{l\in\pi} g^{(|l|)}(t),
\label{eq:resurrexit} % (12)
\end{equation}
where the coefficients $\binom{n}{\la;k}$ of~\eqref{eq:remix} count
partitions yielding equal summands. (Note that
$\sum_{|\la|=k} \binom{n}{\la;k} = \stirlingtwo{n}{k}$ is the Stirling
number of the second kind \cite[Chap.~6]{GrahamKP94}.) For instance,
for $n = 4$ one immediately finds:
\begin{align*}
(f \circ g)^\ivth(t) 
&= f'(g(t))\, g^\ivth(t) + 4 f''(g(t))\, g'(t) g'''(t)
+ 3 f'''(g(t))\, g''(t)^2
\\
&\qquad + 6 f'''(t)\, g'(t)^2 g''(t) + f^\ivth(t)\, g'(t)^4.
\end{align*}
The above can be generalized to formal series in several
indeterminates as follows. Consider the set of maps from
$\{1,\dots,n\}$ to a set of \textit{colours} $\{1',\dots,N'\}$. This
allows consideration of coloured partitions of $\{1,\dots,n\}$, with
monocoloured partitions being of the same type as before. There are
$m$ families of those series, one for each colour, with tangency at
the identity being enforced by:
$$
\del_{t_i} f^i(0,\dots,0) = 1, \quad
\del_{t_j} f^i(0,\dots,0) = 0  \text{ for } j \neq i;
\word{with} i,j \in \{1',\dots,N'\}.
$$
Then the very formula \eqref{eq:resurrexit} is valid, provided we
understand now $|l|$ as a vector of colours. A fully multivariate
treatment in this vein is provided by~\cite{ConstantineS96}. See the
simplest case $N' = 2$ in the foreword to the book~\cite{Primula} --
this book contains much interesting information besides.

The series \eqref{eq:N'-series} can be regarded as characters of
``coloured'' Faà di~Bruno Hopf algebras $\sF(N')$ \cite{HaimanS89}.
For any finite set $X$ gifted with a colouring map
$\theta \: X \to \{1',\dots,N'\}$, one considers partitions~$\pi$
whose sets of blocks are also coloured, provided
$\theta(\{x\}) = \theta(x)$ for singletons. Such coloured partitions
form a poset, with $\pi \leq \rho$ if $\pi$ refines $\rho$ as
partitions, and if $\theta_\pi(B) = \theta_\rho(B)$ for each block $B$
of~$\pi$ which is also a block of~$\rho$; this condition entails that
$\rho$ induces a coloured partition $\rho \stroke \pi$ of the set of
blocks of~$\pi$. Coloured partitions $\pi$ of~$X$ with $\theta(X) = r$
(i.e., the one-block partition of~$X$ is assigned the colour~$r$) form
a poset $\Pi^r_{\bar n}$ where $\bar n \in \bN^{N'}$ counts the
colours of its elements; their types $\tPi^r_{\bar n}$ generate the
Hopf algebra $\sF(N')$, with coproduct given by:
$$
\Dl \tPi^r_{\bar n} := \sum_{\pi\in\Pi^r_{\bar n}} \biggl( 
\prod_{B\in\pi} \tPi^{\theta(B)}_{|B|} \biggr) \ox \tPi^r_{|\pi|}.
$$

A character $f$ of $\sF(N')$ is specified by its values on algebra
generators $f^r_{\bar n} = f(\tPi^r_{\bar n})$, which yield
coefficients of the $N'$-series \eqref{eq:N'-series}. The convolution
of two such characters $g$, $f$ has coefficients
$$
g * f\bigl( \tPi^r_{\bar n} \bigr) = \sum_{\pi\in\Pi^r_{\bar n}} 
f^r_{|\pi|} \prod_{B\in\pi} g^{\theta(B)}_{|B|} 
= \sum_{|\bar k|\leq|\bar n|} f^r_{\bar k}/\bar k!
\prod_{(B_1,\dots,B_{|\bar k|})} 
g^1_{|B_1|}\dots g^{N'}_{|B_{|\bar k|}|}
$$
where the second product ranges over \textit{ordered} coloured
partitions of a set with $|\bar n|$ elements; since there are
$\bar n! \big/ \prod_i (\bar m_i!)^{\la_{i,\bar m_i}}$ of these with
prescribed colours, rearrangement of the right hand side yields a
formula for~\eqref{eq:N'-series}. Thus, the character group of
$\sF(N')$ is anti-isomorphic to the group of $N'$-series
like~\eqref{eq:N'-series} under composition. Also, the antipode on
$\sF(N')$ provides Lagrange reversion in several 
variables~\cite{HaimanS89}.

For the applications of Faà di~Bruno algebras to real analytic
function theory we recommend~\cite{KrantzP02}. The Faà di~Bruno
algebras (perhaps involving functional derivatives) have applications
in quantum field theory. Some elementary ones are described
in~\cite{Quaoar}. Deeper ones related to renormalization theory were
broached in~\cite{CaswellK82,ConnesKr01}, and further explored
in~\cite{Foissy08}. A nice, relatively recent work in this respect
is~\cite{KrYeats17}.

Much remains to be delved into. Noncommutative Bell polynomials were
studied in~\cite{KuruschLM14}. Faà di~Bruno's formulas for operads
have been developed in~\cite{KockW19}. An application to control
theory is found in~\cite{KuruschG17}. A recent comprehensive treatise
on Bell polynomials and generalized Lagrange
inversion~\cite{Schreiber21} is also recommended.

\end{document}